\documentclass[12pt]{amsart}

\usepackage{graphicx}
\usepackage{amsmath}

\usepackage{tikz}

\def\frak{\mathfrak}

\def\C{\mathbb{C}}

\def\R{\mathbb{R}}

\def\cD{\mathcal{D}}
\def\cE{\mathcal{E}}

\def\cP{\mathcal{P}}

\def\al{\alpha}

\def\ga{\gamma}

\def\ka{\kappa}
\def\la{\lambda}

\def\si{\sigma}

\def\om{\omega}

\def\Si{\Sigma}

\def\Om{\Omega}

\newcommand{\der}{{\rm d}}

\numberwithin{equation}{section}

\newtheorem{theorem}{Theorem}[section]

\newtheorem{proposition}[theorem]{Proposition}

\theoremstyle{remark}

\theoremstyle{remark}

\usepackage{amssymb,stmaryrd}
\usepackage{amscd}

\makeatletter
\@namedef{subjclassname@2020}{%
  \textup{2020} Mathematics Subject Classification}
\makeatother

\author{Matthew Randall}

\address{School of Mathematics and Statistics\\
Nanjing University of Information Science and Tecnhology\\
219 Ningliu Road, Nanjing 210044\\
China}
\email{100093@nuist.edu.cn}

\title{Local equivalence of some maximally symmetric rolling distributions and $SU(2)$ Pfaffian systems}


\thanks{The earliest version of the manuscript is supported by the Grant agency of the Czech Republic P201/12/G028.}

\begin{document}

\begin{abstract}
We give a description of Nurowski's conformal structure for some examples of bracket-generating rank 2 distributions in dimension 5, aka $(2,3,5)$-distributions, namely the An-Nurowski circle twistor distribution for pairs of surfaces of constant Gauss curvature rolling without slipping or twisting over each other. In the case of hyperboloid surfaces whose curvature ratios give maximally symmetric $(2,3,5)$-distributions, we find the change of coordinates that map the conformal structure to the flat metric of Engel. We also consider a rank $3$ Pfaffian system in dimension 5 with $SU(2)$ symmetry obtained by rotating two of the 1-forms in the Pfaffian system of the spheres rolling distribution, and discuss complexifications of such distributions.    
\end{abstract}

\maketitle

\pagestyle{myheadings}
\markboth{Randall}{Symmetric rolling distributions and $SU(2)$ Pfaffian systems}

\tableofcontents

\section{Introduction}

Ever since Cartan initiated the study of $(2,3,5)$-distributions over 100 years ago using pioneering methods he developed in \cite{Cartan1910}, which eventually gave rise to the eponymous theory of connections and the method of equivalence, the subject has garnered considerable interest and attention over the passing years. A selection of the vast literature on the development of the subject of $(2,3,5)$-distributions includes \cite{AN14}, \cite{ANp13}, \cite{BH}, \cite{BLN}, \cite{BM}, \cite{BH93-curves}, \cite{B00-dual}, \cite{DH}, \cite{DK}, \cite{GPW}, \cite{HS}, \cite{LNS}, \cite{PN05-conf}, \cite{KS06}, \cite{SW17}, \cite{Strazzullo}, \cite{DTp16}, \cite{TW13} and \cite{TW14}, spurred on by problems involving rolling in mechanics, questions regarding holonomy, parabolic differential geometry and the theory of Lie groups. It also inspired further contributions in the field of subriemannian geometry and control theory \cite{subri-book}. 

A $(2,3,5)$-distribution on a 5 manifold $M^5$ is a maximally non-integrable bracket-generating rank 2 distribution. Such a distribution determines a filtration of the tangent bundle given by
\begin{align*}
\cD \subset [\cD,\cD] \subset [[\cD,\cD],\cD]=TM.
\end{align*} 
The distribution $[\cD,\cD]$ has rank 3 while the full tangent space $TM$ has rank 5, hence such a geometry is also called a $(2,3,5)$-distribution. In Section 5 of \cite{PN05-conf}, it is shown how to associate canonically to such a $(2,3,5)$-distribution a conformal class of metrics of split signature $(2,3)$, henceforth known as Nurowski's conformal structure or Nurowski's conformal metrics, such that the rank 2 distribution is isotropic with respect to any metric in the conformal class. The method of equivalence \cite{Cartan1910} (also see the introduction to \cite{AN14}, Section 5 of \cite{PN05-conf} and \cite{Strazzullo}) applied to $(2,3,5)$-distributions produces the 1-forms $(\theta_1, \theta_2,\theta_3, \theta_4, \theta_5)$ and $(\Om_1, \ldots, \Om_7)$ that together satisfy the structure equations (\ref{cse}) given below. These 1-forms and two additional 1-forms $(\Om_8, \Om_9)$ together define a rank 14 principal bundle over the 5-manifold $M$ (see \cite{Cartan1910} and Section 5 of \cite{PN05-conf}). A representative metric in Nurowski's conformal class \cite{PN05-conf} is given by
\begin{align}\label{metric}
g=2 \theta_1 \theta_5-2\theta_2 \theta_4+\frac{4}{3}\theta_3 \theta_3.
\end{align}
When $g$ has vanishing Weyl tensor, the distribution is called maximally symmetric and has the smallest exceptional Lie group split $G_2$ as its group of local symmetries. For details, see the introduction to \cite{AN14} and Section 5 of \cite{PN05-conf}. Split $G_2$ is also the automorphism group of the split octonions and the relationship between maximally symmetric $(2,3,5)$-distributions and the non-associative division algebra of the split octonions is discussed  in \cite{KS06} and in Section 2 of \cite{TW13}.

$(2,3,5)$-distributions also arise from the study of the configuration space of two surfaces rolling without slipping or twisting over each other \cite{AN14}, \cite{BH} and \cite{BM}. The configuration space can be realised as the An-Nurowski circle twistor distribution \cite{AN14} and in the case of two spheres with radii in the ratio $1:3$ rolling without slipping or twisting over each other, there is again maximal $G_2$ symmetry. 

The motivation of this article is to derive Nurowski's conformal structure(s) for the homogeneous example. Origins of the homogeneous model can be traced to Sections $52$ and $53$ (pages 178-180) of \cite{Cartan1910}, and more recently investigated in the work of \cite{AN14}. In this latter study the 1-forms in the coframe $(\theta_1, \theta_2,\theta_3, \theta_4, \theta_5)$ satisfying the Cartan structure equations (\ref{cse}) are given in Proposition 5 of \cite{AN14}. 

In Section \ref{max1}, we recast in a more simplified form Proposition 5 and write down Nurowski's conformal structure associated to the An-Nurowski circle twistor distribution \cite{AN14} for 1) two surfaces both of positive Gauss curvature in Theorem \ref{thman1} and for 2) two surfaces both of negative Gauss curvature in Theorem \ref{thman2}. 

In Section \ref{max2}, in the case where the surfaces are hyperboloids both of negative Gauss curvature and the curvature ratios are in $1:\frac{1}{9}$ and $1:9$, we find the change of coordinates that bring it to the flat Engel distribution (Propositions \ref{engelco1} and \ref{engelco2}) and write down the vector fields that bracket generate the real split Lie algebra of $\frak{g}_2$ (Theorems \ref{g2asthm} and \ref{g2bsthm}). In these coordinates, the conformal structure is none other than that of the split signature flat metric on $\R^5$ (or $\C^5$) obtained by Engel \cite{Engel} near the end of the 19th century. 

In Section \ref{max3}, we introduce the notion of a Pfaffian system with left-invariant $SU(2)$ symmetry by rotating two of the three 1-forms in the Pfaffian system of the spheres rolling distribution. The vector fields in this rotated distribution becomes $su(2)$-valued. The conformal structure given in terms of this Pfaffian system is presented in Theorems \ref{r1m} and \ref{real3a}.

In Sections \ref{max4} and \ref{max5} we study two different possible complexifications of the $SU(2)$ Pfaffian system and the 1-forms in the rolling distribution. 
In Section \ref{max4}, we study complex Pfaffian systems obtained by allowing the function describing the surface of revolution to be complex. The coframing for Nurowski's metric in this setting gives only constants and therefore defines a homogeneous structure. We obtain the complex analogues of Theorems \ref{r1m} and \ref{real3a} in Section $\ref{max3}$. 

In Section \ref{max5}, we double the $SU(2)$ symmetry of the rotated An-Nurowski circle twistor bundle by considering the Riemannian surface element of the unit sphere arising from one copy of $SU(2)$ and the other Riemannian surface element with Gaussian curvature $9$ or $\frac{1}{9}$ from another copy of $SU(2)$. Both Lie algebras of $su(2)$ are parametrised by the left-invariant vector fields.

We believe the results would be useful to someone who would like to work with the conformal structure associated to maximally symmetric $(2,3,5)$-distributions. 

The computations here are done with the aid of the \texttt{DifferentialGeometry} package in MAPLE 17 and later MAPLE 18.

\section{Nurowski's metric for pairs of surfaces with constant Gauss curvature}\label{max1}

Following \cite{AN14}, let us recall Nurowski's metric associated to the $(2,3,5)$-distribution of a pair of surfaces with constant Gauss curvature rolling without slipping and twisting over each other. In Proposition 5 in pages 407-408 of \cite{AN14}, the conformal structure associated to a pair of surfaces $(\Si^2, \Si'^2)$, one with Killing symmetry and another with constant Gaussian curvature, is provided. Using the coordinates
$x_1=\theta$, $x_2=\phi$, $x_3=x$, $x_4=q$, $\phi=\psi$ and the notation $\chi$ instead of $\si$ in \cite{AN14}, we have
 \begin{align*}
 \chi_1=\der \theta,\qquad \chi_2=f(\theta)\der \phi,\qquad \chi_3=\der x,\qquad \chi_4=h(x)\der q.
 \end{align*}
 The circle twistor distribution is spanned by the vector fields
 \begin{align*}
 V^1&=\partial_\theta+\cos(\psi)\partial_x+\sin(\psi)\left(\frac{1}{h}\partial_q-\frac{h_x}{h}\partial_\psi\right),\\
 V^2&=\frac{1}{f}\partial_\phi+\frac{f_\theta}{f}\partial_\psi-\sin(\psi)\partial_x+\cos(\psi)\left(\frac{1}{h}\partial_q-\frac{h_x}{h}\partial_\psi\right),
 \end{align*}
 which are annihilated by the 1-forms
 \begin{align}\label{an1forms}
 \om_1&=\chi_1-\cos(\psi)\chi_3-\sin(\psi) \chi_4,\nonumber\\
 \om_2&=\chi_2+\sin(\psi)\chi_3-\cos(\psi)\chi_4,\\
  \om_3&=\der\psi+\frac{h_{x}}{h}\chi_4-\frac{f_{\theta}}{f}\chi_2.\nonumber
 \end{align}
 This can be completed to a coframe by taking
 \begin{align*}
 \om_4=\sin(\psi)\chi_3-\cos(\psi) \chi_4 \hspace{12pt}\mbox{and} \hspace{12pt} \om_5=\cos(\psi)\chi_3+\sin(\psi)\chi_4.
 \end{align*}
Based on the results of \cite{AN14}, there are three different scenarios to consider where the $(2,3,5)$-distribution determined by surfaces rolling without slipping and twisting over each other have maximal split $G_2$ symmetry. The first is when both surface elements have positive constant Gauss curvature with ratio of the curvatures $1:\frac{1}{9}$ or $1:9$; the second is when both surface elements have negative constant Gauss curvature with ratio of the curvatures $1:\frac{1}{9}$ or $1:9$; the third is the case of the surfaces of An and Nurowski rolling on a plane.

Here we first give the conformal structure of the $(2,3,5)$-distribution in the cases where: 1) both surfaces are of positive constant Gauss curvature, and 2) both surfaces are of negative constant Gauss curvature. We then specialise to the maximally symmetric cases by taking the above-mentioned ratios. We will leave the discussion of scenario 3) elsewhere. 

When both metrics are of positive constant Gauss curvature, we suppose $f$ and $h$ satisfy the differential equations
 \begin{align*}
 f_{\theta \theta}=-b^2f,\qquad h_{xx}=-c^2h.
 \end{align*}
We assume that the surface element of $\Si^2$
 \begin{align*}
 ds_1^2=\chi_1^2+\chi_2^2=\der \theta^2+f(\theta)^2\der \phi^2
 \end{align*}
 has Gauss curvature $K(ds_1^2)=-\frac{f_{\theta \theta}}{f}=b^2$
 and the surface element of the other surface $\Si'^2$
  \begin{align*}
 ds_2^2=\chi_3^2+\chi_4^2=\der x^2+h(x)^2\der q^2
 \end{align*}
 has Gauss curvature $K(ds_2^2)=-\frac{h_{xx}}{h}=c^2$.

 To obtain the coframing for Nurowski's conformal structure, we take
 \begin{align*}
 \theta_1&=\om_1,\qquad
 \theta_2=\om_2,\qquad
 \theta_3=\frac{1}{(b^2-c^2)^{1/3}}\om_3,\\
 \theta_4&=(b^2-c^2)^{1/3}\om_4-\frac{7b^2-3c^2}{10(b^2-c^2)^{2/3}}\theta_2,\\
 \theta_5&=(b^2-c^2)^{1/3}\om_5+\frac{7b^2-3c^2}{10(b^2-c^2)^{2/3}}\theta_1.
 \end{align*}
 The coframing satisfy Cartan's structure equations \cite{AN14}, \cite{ANp13}, \cite{Cartan1910},  \cite{PN05-conf} given by
\begin{align}\label{cse}
\der \theta_1&=\theta_1\wedge (2\Om_1+\Om_4)+\theta_2\wedge \Om_2+\theta_3 \wedge \theta_4,\nonumber\\
\der \theta_2&=\theta_1\wedge\Om_3+\theta_2\wedge (\Om_1+2\Om_4)+\theta_3 \wedge \theta_5,\nonumber\\
\der \theta_3&=\theta_1\wedge\Om_5+\theta_2\wedge\Om_6+\theta_3\wedge (\Om_1+\Om_4)+\theta_4 \wedge \theta_5,\\
\der \theta_4&=\theta_1\wedge\Om_7+\frac{4}{3}\theta_3\wedge\Om_6+\theta_4\wedge \Om_1+\theta_5 \wedge \Om_2,\nonumber\\
\der \theta_5&=\theta_2\wedge \Om_7-\frac{4}{3}\theta_3\wedge \Om_5+\theta_4\wedge\Om_3+\theta_5\wedge \Om_4.\nonumber
\end{align}

 From this coframing, we obtain
 \begin{align*}
 g&=2\theta_1\theta_5-2\theta_2\theta_4+\frac{4}{3}\theta_3^2\\
 &=2(b^2-c^2)^{1/3}(\om_1\om_5-\om_2\om_4)+\frac{7b^2-3c^2}{5(b^2-c^2)^{2/3}}(\om_1^2+\om_2^2)+\frac{4}{3}\frac{1}{(b^2-c^2)^{2/3}}\om_3^2.
 \end{align*}
 Rescaling by the constant $(b^2-c^2)^{-1/3}$ gives
 \begin{align*}
 g=2(\om_1\om_5-\om_2\om_4)+\frac{7b^2-3c^2}{5(b^2-c^2)}(\om_1\om_1+\om_2\om_2)+\frac{4}{3}\left(\frac{1}{b^2-c^2}\right)\om_3\om_3.
 \end{align*}
 Since we have the following identities
 \begin{align*}
 \om_1&=\chi_1-\om_5,\qquad \om_2=\chi_2+\om_4, \qquad  \om_4^2+\om_5^2=\chi_3^2+\chi_4^2,\\
 \chi_1\om_5&=\frac{1}{2}(\chi_1^2+\om_5^2-(\chi_1-\om_5)^2),\qquad \chi_2\om_4=\frac{1}{2}((\chi_2+\om_4)^2-\chi_2^2-\om_4^2),
 \end{align*}
 we can therefore express the metric more simply as
 \begin{align*}
 g=\chi_1^2+\chi_2^2-(\chi_3^2+\chi_4^2)+\frac{2}{5}\left(\frac{b^2+c^2}{b^2-c^2}\right)(\om_1^2+\om_2^2)+\frac{4}{3}\left(\frac{1}{b^2-c^2}\right)\om_3^2. 
 \end{align*}
 \begin{theorem}[Proposition 5 of \cite{AN14}]\label{thman1}
 Let $\Si^2$ and $\Si'^2$ be two surfaces both of positive constant Gauss curvature, one of curvature $b^2>0$ and the other of curvature $c^2>0$ and with the surface elements parametrised by the functions $f$ satisfying $f_{\theta \theta}=-b^2 f$ and $h$ satisfying $h_{xx}=-c^2 h$ respectively. Nurowski's metric for the An-Nurowski circle twistor distribution of these two surfaces rolling without slipping or twisting over each other is given by
  \begin{align*}
 g=\chi_1^2+\chi_2^2-(\chi_3^2+\chi_4^2)+\frac{2}{5}\left(\frac{b^2+c^2}{b^2-c^2}\right)(\om_1^2+\om_2^2)+\frac{4}{3}\left(\frac{1}{b^2-c^2}\right)\om_3^2. 
 \end{align*}
 The metric is conformally flat and the distribution has maximal split $G_2$ symmetry whenever $b^2-9c^2=0$ or $9b^2-c^2=0$.
 \end{theorem}
 When $b^2=9c^2$, we obtain the metric 
  \begin{align*}
 \tilde g=\chi_1^2+\chi_2^2-(\chi_3^2+\chi_4^2)+\frac{1}{2}(\om_1^2+\om_2^2)+\frac{1}{6c^2}\om_3^2
 \end{align*}
 and when $b^2=\frac{1}{9}c^2$, we have 
  \begin{align*}
 \tilde g=\chi_1^2+\chi_2^2-(\chi_3^2+\chi_4^2)-\frac{1}{2}(\om_1^2+\om_2^2)-\frac{3}{2c^2}\om_3^2.
 \end{align*}
These are the metrics associated to the $(2,3,5)$-distribution of two spheres $(\Si^2, \Si'^2)$ with radii in the ratios of $1:3$ and $3:1$ respectively rolling without slipping and twisting over each other. 

When both surfaces are of negative constant Gauss curvature, we suppose $f$ and $h$ satisfy the differential equations
 \begin{align*}
 f_{\theta \theta}=b^2f,\qquad h_{xx}=c^2h.
 \end{align*}
 It follows that the surface element
 \begin{align*}
 ds_1^2=\der \theta^2+f(\theta)^2\der \phi^2
 \end{align*}
 of $\Si^2$ has Gauss curvature $K(ds_1^2)=-\frac{f_{\theta \theta}}{f}=-b^2$
 and the surface element
  \begin{align*}
 ds_2^2=\der x^2+h(x)^2\der q^2
 \end{align*}
of $\Si'^2$ has Gauss curvature $K(ds_2^2)=-\frac{h_{xx}}{h}=-c^2$.
To obtain the coframing for Nurowski's conformal structure, we take
 \begin{align*}
 \theta_1&=\om_1,\quad \theta_2=\om_2,\quad  \theta_3=-\frac{1}{(b^2-c^2)^{1/3}}\om_3,\\
 \theta_4&=-(b^2-c^2)^{1/3}\om_4+\frac{7b^2-3c^2}{10(b^2-c^2)^{2/3}}\theta_2,\quad \theta_5=-(b^2-c^2)^{1/3}\om_5-\frac{7b^2-3c^2}{10(b^2-c^2)^{2/3}}\theta_1.
 \end{align*}
 From this, we obtain
 \begin{align*}
 g&=2\theta_1\theta_5-2\theta_2\theta_4+\frac{4}{3}\theta_3^2\\
 &=-2(b^2-c^2)^{1/3}(\om_1\om_5-\om_2\om_4)-\frac{7b^2-3c^2}{5(b^2-c^2)^{2/3}}(\om_1^2+\om_2^2)+\frac{4}{3}\frac{1}{(b^2-c^2)^{2/3}}\om_3^2.
 \end{align*}
 Rescaling by the constant $-\frac{1}{(b^2-c^2)^{1/3}}$ gives
 \begin{align*}
 g=2(\om_1\om_5-\om_2\om_4)+\frac{7b^2-3c^2}{5(b^2-c^2)}(\om_1\om_1+\om_2\om_2)-\frac{4}{3}\left(\frac{1}{b^2-c^2}\right)\om_3\om_3.
 \end{align*}
Again using the above identities we can express the metric more simply as
 \begin{align*}
 g=\chi_1^2+\chi_2^2-(\chi_3^2+\chi_4^2)+\frac{2}{5}\left(\frac{b^2+c^2}{b^2-c^2}\right)(\om_1^2+\om_2^2)-\frac{4}{3}\left(\frac{1}{b^2-c^2}\right)\om_3^2. 
 \end{align*}
 \begin{theorem}[Proposition 5 of \cite{AN14}]\label{thman2}
  Let $\Si^2$ and $\Si'^2$ be two surfaces both of negative constant Gauss curvature, one of curvature $b^2<0$ and the other of curvature $c^2<0$ and with the surface elements parametrised by the functions $f$ satisfying $f_{\theta \theta}=b^2 f$ and $h$ satisfying $h_{xx}=c^2 h$ respectively. 
 Nurowski's metric for the An-Nurowski circle twistor distribution of these two surfaces rolling without slipping or twisting over each other is given by
  \begin{align*}
 g=\chi_1^2+\chi_2^2-(\chi_3^2+\chi_4^2)+\frac{2}{5}\left(\frac{b^2+c^2}{b^2-c^2}\right)(\om_1^2+\om_2^2)-\frac{4}{3}\left(\frac{1}{b^2-c^2}\right)\om_3^2. 
 \end{align*}
 The metric is conformally flat and the distribution has maximal split $G_2$ symmetry whenever $b^2-9c^2=0$ or $9b^2-c^2=0$.
 \end{theorem}
  When $b^2=9c^2$, we obtain the metric 
  \begin{align}\label{ga-hyperboloid}
 \tilde g=\chi_1^2+\chi_2^2-(\chi_3^2+\chi_4^2)+\frac{1}{2}(\om_1^2+\om_2^2)-\frac{1}{6c^2}\om_3^2
 \end{align}
 and when $b^2=\frac{1}{9}c^2$, we have 
  \begin{align}\label{gb-hyperboloid}
 \tilde g=\chi_1^2+\chi_2^2-(\chi_3^2+\chi_4^2)-\frac{1}{2}(\om_1^2+\om_2^2)+\frac{3}{2c^2}\om_3^2.
 \end{align}
The metric becomes singular when the ratios of the radii and curvatures are $1:1$. In this case, the distribution is integrable and no longer a $(2,3,5)$-distribution.

 \section{Maximally symmetric twistor distribution for surfaces with negative constant Gauss curvature}\label{max2}
 In this section, we consider the $(2,3,5)$-distribution \`a la An and Nurowski \cite{AN14} in the case of hyperboloids of constant negative Gauss curvature rolling without slipping or twisting over each other. In the case where maximal symmetry is obtained when the ratio of the curvatures is either $1:\frac{1}{9}$ or $1:9$, we find the change of coordinates to the flat Engel distribution and write down the Lie algebra of vector fields that bracket-generate the split real form of $\frak{g}_2$. In order to accomplish this, we specialise to the case where $f(\theta)=e^{\theta}$ so that $b=1$ and $h(x)=e^{c x}$ where $c^2=\frac{1}{9}$ or $9$. We consider firstly the case where $c^2=\frac{1}{9}$ by taking $c=\frac{1}{3}$. The situation where $c=-\frac{1}{3}$ is similar. 
 
 The surface metric given by
 \begin{align*}
 ds_1^2=\der \theta^2+e^{2 \theta} \der \phi^2
 \end{align*}
 has constant negative Gauss curvature $-1$. 
 The surface metric given by
 \begin{align*}
 ds_2^2=\der x^2+e^{\frac{2x}{3}} \der q^2
 \end{align*}
 has constant negative Gauss curvature $-\frac{1}{9}$. 
 
The distribution is spanned by the vector fields
\begin{align*}
V^1&=\partial_\theta+\cos(\psi)\partial_x+\sin(\psi)(e^{-\frac{x}{3}}\partial_{q}-\frac{1}{3}\partial_{\psi}),\\
V^2&=e^{-\theta}\partial_\phi+\partial_\psi-\sin(\psi)\partial_{x}+\cos(\psi)(e^{-\frac{x}{3}}\partial_{q}-\frac{1}{3}\partial_{\psi}),
 \end{align*} 
 which are annihilated by the 1-forms
  \begin{align*}
 \om_1&=\der \theta-\cos(\psi) \der x-\sin(\psi)e^{\frac{x}{3}}\der q,\\
 \om_2&=e^{\theta} \der \phi+\sin(\psi)\der x-\cos(\psi)  e^{\frac{x}{3}}\der q,\\ 
 \om_3&=\der \psi-e^\theta \der \phi+\frac{1}{3}e^\frac{x}{3}\der q.
 \end{align*}
 
 To write down the real vector fields that bracket-generate split $\frak{g}_2$ and relate the maximally symmetric rolling $(2,3,5)$-distribution to the flat Engel distribution, we take the coordinate transformation given by
 \begin{align*}
c_1&=4\cos(\frac{\psi}{2})e^{\frac{\theta}{2}-\frac{x}{2}}-\frac{16}{3}\sin(\frac{\psi}{2})e^{\frac{\theta}{2}-\frac{x}{6}}q,\\
c_2&=-\frac{8}{3}qe^{-\frac{x}{6}}(e^{\frac{\theta}{2}}\sin(\frac{\psi}{2})\phi+e^{-\frac{\theta}{2}}\cos(\frac{\psi}{2}))+2e^{-\frac{x}{2}}(e^{\frac{\theta}{2}}\cos(\frac{\psi}{2})\phi-e^{-\frac{\theta}{2}}\sin(\frac{\psi}{2})),\\
c_3&=\frac{2}{3}q,\\
c_4&=-2e^{\frac{\theta}{2}-\frac{x}{6}}\sin(\frac{\psi}{2}),\\
c_5&=e^{-\frac{x}{6}}(e^{-\frac{\theta}{2}}\cos(\frac{\psi}{2})+e^{\frac{\theta}{2}}\sin(\frac{\psi}{2})\phi).
 \end{align*}
 It follows that the 1-forms $\Theta_1$, $\Theta_2$, $\Theta_3$ given by
 \begin{align*}
\Theta_1&:=\der c_1-2c_4 \der c_3-4c_3 \der c_4\\
&=2e^{\frac{\theta}{2}-\frac{x}{2}}(\cos(\frac{\psi}{2})\om_1-\sin(\frac{\psi}{2})(\om_2+\om_3)),\\
\Theta_2&:=\der c_2+2c_5 \der c_3+4c_3 \der c_5\\
&=e^{-\frac{x}{2}}(e^{\frac{\theta}{2}}\phi(\cos(\frac{\psi}{2})\om_1-\sin(\frac{\psi}{2})(\om_2+\om_3))+e^{-\frac{\theta}{2}}(\cos(\frac{\psi}{2})(\om_2-\om_3)+\sin(\frac{\psi}{2})\om_1)),\\
\Theta_3&:=\der c_3+c_5 \der c_4-c_4 \der c_5\\
&=-e^{-\frac{x}{3}}(\om_3+\sin(\psi)\om_1+\cos(\psi)\om_2)
 \end{align*} 
 are in the span of the Pfaffian system $\{\om_1,\om_2,\om_3\}$.

Following \cite{R21b}, we write down the vector fields
\begin{align*}
Z^1=\partial_{c_3}+2c_5\partial_{c_2}-2c_4\partial_{c_1},\\
Z^2=\partial_{c_4}+4c_3\partial_{c_1}-2c_5\partial_{c_3},\\
Z^3=\partial_{c_5}+2c_4\partial_{c_3}-4c_3\partial_{c_2},
\end{align*}
and define
\begin{align*}
S^1:=Z^2+c_5Z^1,\hspace{12pt} S^2:=Z^3-c_4Z^1, \hspace{12pt} S^3:=-c_1Z^2+c_2Z^3-(c_1c_5+c_2c_4+c_3^2)Z^1.
\end{align*}
It follows from the formulas of $(c_1, c_2, c_3, c_4, c_5)$ given above that 
\begin{align}\label{g2a1}
S^1&=\frac{3}{4}\sin(\frac{\psi}{2})e^{\frac{x}{6}-\frac{\theta}{2}}X^1-\frac{3}{2}\cos(\frac{\psi}{2})e^{\frac{x}{6}-\frac{\theta}{2}}X^2+\frac{1}{2}\phi S^2,\\ \label{g2a2}
S^2&=-\frac{3}{2}\cos(\frac{\psi}{2})e^{\frac{\theta}{2}+\frac{x}{6}}X^1-3\sin(\frac{\psi}{2})e^{\frac{\theta}{2}+\frac{x}{6}}X^2,\\ \label{g2a3}
S^3&=4qX^1+6e^{-\frac{x}{3}}X^2-\frac{4}{9}q^2 Z^1,
\end{align}
where
\begin{align*}
X^1=\cos(\psi)V^1-\sin(\psi)V^2,\qquad X^2=\cos(\psi) V^2+\sin(\psi) V^1
\end{align*}
are the vector fields obtained from rotating the distribution
and
\begin{align*}
Z^1&=\frac{3}{2}\partial_q-3e^{\frac{x}{3}}\partial_{\psi}+3e^{\frac{x}{3}}(\sin(\psi)\partial_\theta+\cos(\psi)(e^{-\theta}\partial_{\phi}+\partial_{\psi})).
\end{align*}
The vector fields $S^1$ and $S^2$ are in the span of the distribution and are annihilated by the 1-forms $\{\Theta_1, \Theta_2, \Theta_3\}$.

The vector fields $S^1$, $S^2$ and $S^3$ given by (\ref{g2a1}), (\ref{g2a2}) and (\ref{g2a3}) pairwise bracket-generate the Lie algebra of split ${\frak g}_2$ in the following sense: 
defining the following vector fields by the Lie brackets
\begin{align*}
S^4&:=[S^1,S^2],\quad S^5:=[S^2,S^3],\quad S^6:=[S^3,S^1],\\
L^1&:=[S^1,S^4],\quad L^3:=[S^2,S^5],\quad L^5:=[S^3,S^6],\\
L^2&:=[S^4,S^2],\quad L^4:=[S^5,S^3],\quad L^6:=[S^6,S^1]
\end{align*}
and taking
\begin{align*}
H:=[S^2,S^6], \quad h:=[S^4,S^3],
\end{align*}
we require that the set of vector fields
\begin{align*}
\{S^1,S^2,S^3,S^4,S^5,S^6,\frac{1}{4}(h-H),\frac{\sqrt{3}}{12}(h+H),L^1,L^2,L^3,L^4,L^5,L^6\}
\end{align*}
form the 14-dimensional Lie algebra of split ${\frak g}_2$. With this choice of the Cartan subalgebra spanned by $\frac{1}{4}(h-H)$ and $\frac{\sqrt{3}}{12}(h+H)$, the root diagram is given by the picture below. 
\begin{figure}[h!]
\begin{tikzpicture}
	\draw [stealth-stealth](-1,0) -- (1,0);
\draw (1,0) node[anchor=west] {{\tiny $S^1$}};
\draw (-1,0) node[anchor=east] {{\tiny $S^5$}};
	\draw [stealth-stealth](0,-1.732) -- (0,1.732);
\draw (0,-1.732) node[anchor=north] {{\tiny $L^5$}};
\draw (0,1.732) node[anchor=south] {{\tiny $L^2$}};
\draw [stealth-stealth](-0.5,-0.866) -- (0.5,0.866);
\draw (-0.5,-0.866) node[anchor=north] {{\tiny $S^{3}$}};
\draw (0.5,0.866) node[anchor=south] {{\tiny $S^4$}};
\draw [stealth-stealth](-1.5,-0.866) -- (1.5,0.866);
\draw (-1.5,-0.866) node[anchor=north] {{\tiny $L^4$}};
\draw (1.5,0.866) node[anchor=south] {{\tiny $L^1$}};
\draw [stealth-stealth](1.5,-0.866) -- (-1.5,0.866);
\draw(1.5,-0.866) node[anchor=north] {{\tiny $L^6$}};
\draw (-1.5,0.866) node[anchor=south] {{\tiny $L^3$}};
\draw [stealth-stealth](0.5,-0.866) -- (-0.5,0.866);
\draw (0.5,-0.866) node[anchor=north] {{\tiny $S^{6}$}};
\draw (-0.5,0.866) node[anchor=south] {{\tiny $S^2$}};
\end{tikzpicture}
\end{figure} 

\begin{theorem}\label{g2asthm}
The vector fields $S^1$, $S^2$ and $S^3$ given by (\ref{g2a1}), (\ref{g2a2}), (\ref{g2a3}) pairwise bracket-generate the Lie algebra of split $\frak{g}_2$.
 \end{theorem}    
  The rest of the vector fields can be determined from the corresponding formulas of $(c_1,c_2,c_3,c_4,c_5)$. For instance, we find 
\begin{align*}
 H=6(\partial_x-\frac{1}{3}q\partial_q-\phi\partial_\phi+\partial_\theta)\hspace{12pt} \mbox{and} \hspace{12pt}
 h=4(3\partial_x-q\partial_q).
\end{align*}      

The coordinates for the flat metric are obtained by passing to the flat Engel distribution by taking 
$\frak{r_1}=c_5$,$\frak{r_2}=c_4$, $\frak{r_3}=c_3$, $\frak{r_4}=\frac{1}{2}(c_2+3c_3c_5)$, $\frak{r_5}=\frac{1}{2}(c_1-3c_3c_4)$.
It follows that Nurowski's metric in this instance is conformally related to none other than the flat metric
\begin{align*}
ds^2=2\der \frak{r_1}\der \frak{r_5}+2\der \frak{r_2}\der \frak{r_4}+\der \frak{r_3}^2
\end{align*}
on $\R^5$ (or $\C^5$) determined by Engel in \cite{Engel} and \cite{Engel2} over a century ago and the 1-forms 
\begin{align*}
\der\frak{r_3}+\frak{r_1}\der \frak{r_2}-\frak{r_2}\der \frak{r_1}&=\Theta_3,\\
\der\frak{r_4}+\frac{1}{2}(\frak{r_3}\der \frak{r_1}-\frak{r_1}\der \frak{r_3})&=\frac{1}{2}\Theta_2,\\
\der\frak{r_5}+\frac{1}{2}(\frak{r_2}\der \frak{r_3}-\frak{r_3}\der \frak{r_2})&=\frac{1}{2}\Theta_1
\end{align*}
obtained in \cite{Engel} and \cite{Engel2} are those annihilating the flat Engel distribution spanned by the vector fields
\begin{align*}
T^2-\frac{1}{2}\frak{r}_1T^3\hspace{12pt} \mbox{and} \hspace{12pt} T^1+\frac{1}{2}\frak{r}_2T^3
\end{align*}
where 
\begin{align*}
T^1&=\partial_{\frak{r}_1}+\frac{1}{2}\frak{r}_2\partial_{\frak{r}_3}-\frac{1}{2}\frak{r}_3\partial_{\frak{r}_4},\quad  T^2=\partial_{\frak{r}_2}+\frac{1}{2}\frak{r}_3\partial_{\frak{r}_5}-\frac{1}{2}\frak{r}_1\partial_{\frak{r}_3},\quad T^3=\partial_{\frak{r}_3}+\frak{r}_1\partial_{\frak{r}_4}-\frak{r}_2\partial_{\frak{r}_5}.
\end{align*}

 \begin{proposition}\label{engelco1}
 The coordinate functions 
 \begin{align*}
 \frak{r}_1&=e^{-\frac{x}{6}}(e^{-\frac{\theta}{2}}\cos(\frac{\psi}{2})+e^{\frac{\theta}{2}}\sin(\frac{\psi}{2})\phi),\\
 \frak{r}_2&=-2e^{\frac{\theta}{2}-\frac{x}{6}}\sin(\frac{\psi}{2}),\\
 \frak{r}_3&=\frac{2}{3}q,\\
 \frak{r}_4&=\frac{1}{2}(c_2+3c_3c_5)=-\frac{1}{3}qe^{-\frac{x}{6}}(e^{\frac{\theta}{2}}\sin(\frac{\psi}{2})\phi+e^{-\frac{\theta}{2}}\cos(\frac{\psi}{2}))+e^{-\frac{x}{2}}(e^{\frac{\theta}{2}}\cos(\frac{\psi}{2})\phi-e^{-\frac{\theta}{2}}\sin(\frac{\psi}{2})),\\
 \frak{r}_5&=\frac{1}{2}(c_1-3c_3c_4) =2\cos(\frac{\psi}{2})e^{\frac{\theta}{2}-\frac{x}{2}}-\frac{2}{3}\sin(\frac{\psi}{2})e^{\frac{\theta}{2}-\frac{x}{6}}q
  \end{align*}
    map the Nurowski conformal structure of the maximally symmetric rolling distribution of two hyperboloid surfaces with curvatures in the ratio $1:\frac{1}{9}$ to the flat metric given by Engel. 
    We have
    \begin{align*}
    ds^2=-\frac{2}{3}e^{-\frac{2}{3}x}\tilde g
    \end{align*}
    where $\tilde g$ is given in (\ref{ga-hyperboloid}) with $c=\frac{1}{3}$.
  \end{proposition}  
 
In the other case where $9b^2=c^2$, taking $b=1$ and $c=3$, we have the transformation given by
 \begin{align*}
c_1&=\frac{1}{3}\cos(\frac{\psi}{2})e^{-\frac{3\theta}{2}+\frac{3x}{2}}+\frac{4}{3}\sin(\frac{\psi}{2})e^{-\frac{\theta}{2}+\frac{3x}{2}}\phi,\\
c_2&=\frac{1}{3}(e^{-\frac{\theta}{2}}(\frac{4}{3}e^{-\frac{3x}{2}}\cos(\frac{\psi}{2})-4e^{\frac{3x}{2}}\sin(\frac{\psi}{2}) q)\phi-e^{-\frac{3\theta}{2}}(\frac{1}{3}e^{-\frac{3x}{2}}\sin(\frac{\psi}{2})+e^{\frac{3x}{2}}\cos(\frac{\psi}{2}) q)),\\
c_3&=-\frac{1}{3}\phi,\\
c_4&=-e^{\frac{3x}{2}-\frac{\theta}{2}}\sin(\frac{\psi}{2}),\\
c_5&=e^{-\frac{3x}{2}-\frac{\theta}{2}}(\frac{1}{3}\cos(\frac{\psi}{2})-\sin(\frac{\psi}{2})q e^{3x}).
 \end{align*}
 It follows that the 1-forms $\Theta_1$, $\Theta_2$, $\Theta_3$ are in the span of the 1-forms
  \begin{align*}
 \om_1&=\der \theta-\cos(\psi) \der x-\sin(\psi)e^{3x}\der q,\\
 \om_2&=e^{\theta} \der \phi+\sin(\psi)\der x-\cos(\psi)  e^{3x}\der q,\\ 
 \om_3&=\der \psi-e^\theta \der \phi+3e^{3x}\der q.
 \end{align*}
 We find
 \begin{align*}
 \Theta_1&=-\frac{1}{2}e^{\frac{3}{2}(x-\theta)}(\cos(\frac{\psi}{2})\om_1-\sin(\frac{\psi}{2})(\om_2-\frac{1}{3}\om_3)),\\
 \Theta_2&=\frac{1}{2}\cos(\frac{\psi}{2})e^{-\frac{3\theta}{2}}\left(qe^{\frac{3x}{2}}\om_1+\frac{1}{3}e^{-\frac{3x}{2}}\om_2-\frac{1}{9}e^{-\frac{3x}{2}}\om_3\right)\\
 &\quad{}+\frac{1}{2}\sin(\frac{\psi}{2})e^{-\frac{3\theta}{2}}\left(\frac{1}{3}e^{-\frac{3x}{2}}\om_1-qe^{\frac{3x}{2}}\om_2+\frac{1}{3}qe^{\frac{3x}{2}}\om_3\right),\\
 \Theta_3&=-\frac{1}{2}e^{-\theta}(\om_2+\frac{1}{3}\om_3).
 \end{align*}
 The surface metric given by
 \begin{align*}
 ds_2^2=\der x^2+e^{6x} \der q^2
 \end{align*}
 has constant negative Gauss curvature $-9$.
The distribution is spanned by the vector fields
\begin{align*}
V^1&=\partial_\theta+\cos(\psi)\partial_x+\sin(\psi)(e^{-3x}\partial_q-3\partial_{\psi}),\\
V^2&=e^{-\theta}\partial_\phi+\partial_\psi-\sin(\psi)\partial_x+\cos(\psi)(e^{-3x}\partial_q-3\partial_{\psi}).
 \end{align*} 

Again using the coordinate functions, we can deduce that the vector fields
\begin{align}\label{g2b1}
S^1&=\frac{1}{2}e^{\frac{\theta}{2}-\frac{3x}{2}}\sin(\frac{\psi}{2})V^1+e^{\frac{\theta}{2}-\frac{3x}{2}}\cos(\frac{\psi}{2})V^2-qS^2,\\ \label{g2b2}
S^2&=-\frac{3}{2}e^{\frac{\theta}{2}+\frac{3x}{2}}\cos(\frac{\psi}{2})V^1+3e^{\frac{\theta}{2}+\frac{3x}{2}}\sin(\frac{\psi}{2})V^2,\\ \label{g2b3}
S^3&=-\frac{2}{3}\phi V^1-\frac{1}{3}e^{-\theta}V^2-\frac{1}{9}\phi^2 Z^1,
\end{align}
where
\begin{align*}
Z^1&=-3\partial_\phi-18e^{\theta}\partial_{\psi}-6e^{\theta}(\cos(\psi)(e^{-3x}\partial_q-3 \partial_{\psi})-\sin(\psi)\partial_x)
\end{align*}
pairwise bracket-generate the Lie algebra of split $\frak{g}_2$.

 \begin{theorem}\label{g2bsthm}
The vector fields $S^1$, $S^2$ and $S^3$ given by (\ref{g2b1}), (\ref{g2b2}) and (\ref{g2b3}) pairwise bracket-generate the Lie algebra of split $\frak{g}_2$.
 \end{theorem} 
 Again, the coordinate functions also map the conformal structure to the flat metric of the Engel distribution.
 
   \begin{proposition}\label{engelco2}
 The coordinate functions 
 \begin{align*}
 \frak{r}_1&=-e^{\frac{3x}{2}-\frac{\theta}{2}}\sin(\frac{\psi}{2})q+\frac{1}{3}e^{-\frac{3x}{2}-\frac{\theta}{2}}\cos(\frac{\psi}{2}),\\
 \frak{r}_2&=-e^{\frac{3x}{2}-\frac{\theta}{2}}\sin(\frac{\psi}{2}),\\
 \frak{r}_3&=-\frac{1}{3}\phi,\\
 \frak{r}_4&=(-\frac{1}{6}e^{\frac{3x}{2}-\frac{\theta}{2}}\sin(\frac{\psi}{2})q+\frac{1}{18}e^{-\frac{3x}{2}-\frac{\theta}{2}}\cos(\frac{\psi}{2}))\phi-\frac{1}{6}e^{\frac{3x}{2}-\frac{3\theta}{2}}\cos(\frac{\psi}{2})q-\frac{1}{18}e^{-\frac{3x}{2}-\frac{3\theta}{2}}\sin(\frac{\psi}{2}),\\
 \frak{r}_5&=\frac{1}{6}e^{\frac{3x}{2}}(e^{-\frac{3\theta}{2}}\cos(\frac{\psi}{2})+\phi e^{-\frac{\theta}{2}}\sin(\frac{\psi}{2}))
  \end{align*}
      map the Nurowski conformal structure of the maximally symmetric rolling distribution of two hyperboloid surfaces in the ratio $1:9$ to the flat metric given by Engel. 
    We have
    \begin{align*}
    ds^2=\frac{1}{6}e^{-2\theta}\tilde g
    \end{align*}
    where $\tilde g$ is given in (\ref{gb-hyperboloid}) with $c=3$.
  \end{proposition}
  
In the work of \cite{R21a} Pfaffian systems with real $SL(2)$ symmetry were studied, with the caveat that the metrics on the hyperboloid surfaces are of split signature. By using a parametrisation of the Lie algebra of $sl_2$ arising from the second prolongation of the group action of unimodular fractional linear transformations, a Monge normal form for the rolling distribution with hyperboloid surfaces was also obtained. We now investigate the case where we have real $SU(2)$ symmetry.
  
\section{Spheres rolling and real $SU(2)$ Pfaffian systems}\label{max3}
In this section we consider the case of rolling spheres where we take $f(\theta)=-\sin(\theta)$ and $h(x)=\cos(c x)$ with either $c=\frac{1}{3}$ or $c=3$. In these cases, the distribution has maximal symmetry but the author is unable to find the change of coordinates to the flat Engel distribution. However, we can relate it to real Pfaffian systems with $SU(2)$ symmetry in the following manner.

Let $\si_1$, $\si_2$, $\si_3$ be the left-invariant 1-forms on $SU(2)$, which satisfy the relations
\begin{align*}
\der \si_1=\si_2 \wedge \si_3,\\
\der \si_2=\si_3 \wedge \si_1,\\
\der \si_3=\si_1 \wedge \si_2.
\end{align*}
As a manifold, $SU(2)$ is isomorphic to $S^3$ and these left-invariant 1-forms define a coframing on $S^3$ that can be parameterised by Euler angles $(\psi,\theta,\phi)$ as follows:
\begin{align*}
\si_1&=\sin(\psi)\der \theta-\cos(\psi)\sin(\theta)\der \phi,\\
\si_2&=\cos(\psi)\der \theta+\sin(\psi)\sin(\theta)\der \phi,\\
\si_3&=-\der \psi-\cos(\theta)\der \phi.
\end{align*}
They are dual to the left-invariant vector fields on $S^3$, which we denote by $E^i$ for $i=1,2$ or $3$. In coordinates these vector fields are given by
\begin{align*}
E^1&=\sin(\psi)\partial_\theta+\cos(\psi)\left(\cot(\theta)\partial_{\psi}-\frac{1}{\sin(\theta)}\partial_\phi\right),\\
E^2&=\cos(\psi)\partial_\theta-\sin(\psi)\left(\cot(\theta)\partial_{\psi}-\frac{1}{\sin(\theta)}\partial_\phi\right),\\
E^3&=-\partial_{\psi}.
\end{align*}  
 
 We see that given the 1-forms in (\ref{an1forms}) with $f(\theta)=-\sin(\theta)$, we have
 \begin{align*}
 \sin(\psi)\om_1+\cos(\psi) \om_2&=\sin(\psi) \der \theta-\cos(\psi) \sin(\theta) \der \phi-h \der q=\si_1-h \der q,\\
 \cos(\psi) \om_1-\sin(\psi)\om_2&=\cos(\psi) \der \theta+\sin(\psi) \sin(\theta) \der \phi-\der x=\si_2-\der x,\\
 -\om_3&=-\der \psi-\cos(\theta)\der \phi-h_x\der q=\si_3-h_x \der q.
 \end{align*}
 In other words the $SU(2)$ Pfaffian system is obtained by rotating the 1-forms $\om_1$ and $\om_2$ from the annihilator of the An-Nurowski circle twistor distribution (\ref{an1forms}) with $f(\theta)=-\sin(\theta)$. The vector fields that span the distribution now becomes $su(2)$-valued. 
We shall denote by $o_1$, $o_2$, $o_3$ the 1-forms in the $SU(2)$ Pfaffian system, given by
 \begin{align*}
o_1=\si_1-h(x)\der q,\qquad o_2=\si_2-\der x,\qquad o_3=\si_3-h_x\der q.
\end{align*}
The parametrisation taking $f(\theta)=-\sin(\theta)$ and $h(x)=e^{i c x}$ where $c=\pm \frac{1}{3}$ or $\pm 3$ will be discussed in the next section. In the remaining part of this paper, we treat the case of real left-invariant vector fields on $SU(2)$, without regards to finding the change of coordinates to the flat Engel distribution, which the author is unable to do.  Another generalisation where we express the coframing on $\Si'^2$ in terms of the basis of left-invariant 1-forms for another copy of $S^3$, is given in Section \ref{max5}. 
 
By introducing the notion of a ``sign-reversed" distribution $\tilde \cD$ to be defined in this section, we can rewrite the metrics obtained in Theorem \ref{thman1} in a more symmetric form, namely as a diagonal product of the metric on the sphere $\Si'^2$ and a symmetric bilinear form on $TM/\cD$ and $TM/\tilde \cD$ using the 1-forms in the Pfaffian system that span the complementary bundle $\cD^{\perp}=(TM/\cD)^*$ and $\tilde \cD^{\perp}=(TM/\tilde \cD)^*$.

 In the case where $b^2=1$, $c^2=\frac{1}{9}$, we take $h(x)=\cos(\frac{x}{3})$.
Consider the real Pfaffian system $\cP_{\R}$ given by
the 1-forms
\begin{align*}
o_1=\si_1-\cos\left(\frac{x}{3}\right)\der q,\qquad o_2=\si_2-\der x,\qquad o_3=\si_3+\frac{1}{3}\sin\left(\frac{x}{3}\right)\der q,
\end{align*}
with $o_4=-\chi_3=-\der x$ and $o_5=\chi_4=\cos(\frac{x}{3})\der q$.
This Pfaffian system annihilates the vector fields 
\begin{align*}
X^1=\frac{\partial}{\partial x}+E^2,\quad X^2=\frac{\partial}{\partial q}+\cos\left(\frac{x}{3}\right)E^1-\frac{1}{3}\sin\left(\frac{x}{3}\right)E^3,
\end{align*}
which span the real distribution $\cD_{\R}$.
Observe that
\begin{align*}
o_1^2+o_2^2=\om_1^2+\om_2^2,\qquad o_3=-\om_3,\qquad o_4^2+o_5^2=\om_4^2+\om_5^2=\chi_3^2+\chi_4^2.
\end{align*}

This gives the conformal structure with metric representative
  \begin{align*}
 g_{\R}&=\chi_1^2+\chi_2^2-(\chi_3^2+\chi_4^2)+\frac{1}{2}(\om_1^2+\om_2^2)+\frac{3}{2}\om_3^2\\
 &=\si_1^2+\si_2^2-(o_4^2+o_5^2)+\frac{1}{2}(o_1^2+o_2^2)+\frac{3}{2}o_3^2, 
 \end{align*}
which is conformally flat.

Instead of writing the metric as a product on spheres, we can define a ``sign-reversed'' Pfaffian system $\tilde \cP_{\R}$ corresponding to the symmetry $(x,q,c) \mapsto (-x,-q,-c)$ where $c=\frac{1}{3}$ by reversing signs in $o_1$ and $o_2$, so that
\begin{align*}
\tilde o_1=\si_1+\cos\left(\frac{x}{3}\right)\der q,\qquad
\tilde o_2=\si_2+\der x,\qquad
\tilde o_3=o_3
\end{align*}
and taking $\tilde o_4=\der x=-o_4$, $\tilde o_5=-\cos(\frac{x}{3}) \der q=-o_5$.

The ``sign-reversed" Pfaffian system $\tilde \cP_{\R}$ annihilates the vector fields given by
\begin{align*}
\tilde X^1=\frac{\partial}{\partial x}-E^2,\qquad
\tilde X^2=\frac{\partial}{\partial q}-\cos\left(\frac{x}{3}\right)E^1-\frac{1}{3}\sin\left(\frac{x}{3}\right)E^3,
\end{align*}
which span the ``sign-reversed" real distribution $\tilde \cD_{\R}$.

The ``sign-reversed'' 1-forms give the polarisation identity
\[
\frac{1}{2}(o_1o_1+o_2o_2+\tilde o_1\tilde o_1+\tilde o_2\tilde o_2)=\si_1^2+\si_2^2+o_4o_4+o_5o_5
\]
describing the sum of two sphere elements, with 
\[
\si_1^2+\si_2^2=\chi_1^2+\chi_2^2=\der \theta^2+\sin(\theta)^2\der \phi^2
\] 
the metric of the unit sphere $\Si^2$ and 
\begin{align*}
o_4o_4+o_5o_5=\der x^2+\cos\left(\frac{x}{3}\right)^2\der q^2
\end{align*}
the metric on $\Si'^2$ the sphere of radius $3$.

We can express the above metric now simply as
\begin{align*}
g_{\R}&=-2(o_4 o_4+o_5o_5)+\frac{1}{2}(\tilde{o}_1\tilde{o}_1+\tilde{o}_2\tilde{o}_2+\tilde{o}_3\tilde{o}_3)+o_1o_1+o_2o_2+o_3o_3.
\end{align*}
Observe that this is a product of the metric $o_4 o_4+o_5 o_5$ defined on the surface $\Si'^2$ and the bilinear form 
\begin{align*}
\frac{1}{2}(\tilde{o}_1\tilde{o}_1+\tilde{o}_2\tilde{o}_2+\tilde o_3\tilde o_3)+o_1o_1+o_2o_2+o_3o_3
\end{align*}
on $(TM/\tilde \cD) \times (TM/\cD)$.
The metric tensor associated to $\cD^{\perp}$ is defined by 
\[
g_{\cD^{\perp}}:TM/\cD \times TM/ \cD \rightarrow \R 
\]
and given in the dual basis by
\[
g_{\cD^{\perp}}=o_1o_1+o_2o_2+o_3 o_3,
\]
while the tensor
$g_{\tilde \cD ^{\perp}}:TM/\tilde \cD \times TM/\tilde \cD \rightarrow \R$ is given in the dual basis by
\[
g_{\tilde \cD^{\perp}}=\frac{1}{2}(\tilde{o}_1\tilde{o}_1+\tilde{o}_2\tilde{o}_2+\tilde o_3\tilde o_3).
\]

We have the following result.
\begin{theorem}\label{r1m}
A representative metric of Nurowski's conformal class associated to the maximally symmetric real distribution $\cD_{\R}$ given by the span of 
\begin{align*}
X^1=\frac{\partial}{\partial x}+E^2,\qquad X^2=\frac{\partial}{\partial q}+\cos\left(\frac{x}{3}\right)E^1-\frac{1}{3}\sin\left(\frac{x}{3}\right)E^3,
\end{align*}
is given by 
\begin{align}\label{m1f}
g_{\R}&=-2(o_4 o_4+o_5o_5)+\frac{1}{2}(\tilde{o}_1\tilde{o}_1+\tilde{o}_2\tilde{o}_2+\tilde o_3\tilde o_3)+ (o_1o_1+o_2 o_2+o_3 o_3)\\
&=\chi_1^2+\chi_2^2-(\chi_3^2+\chi_4^2)+\frac{1}{2}(\om_1^2+\om_2^2)+\frac{3}{2}\om_3^2\nonumber
\end{align}
where
$o_1$, $o_2$, $o_3$ are the 1-forms in the Pfaffian system $\cP_{\R}$, which span $(\cD_{\R})^{\perp}$ and $\tilde o_1$, $\tilde o_2$, $\tilde o_3$ are 1-forms in the Pfaffian system $\tilde \cP_{\R}$, which span $(\tilde \cD_{\R})^{\perp}$. The conformal structure determined by $g_{\R}$ is conformally flat and the real distribution has split real form of $G_2$ as its group of symmetries. 
\end{theorem}

For the metric obtained in Theorem \ref{r1m}, we are able to find an Einstein representative in the conformal class. It is given by the following. 
\begin{proposition}
Let $\la$ be a non-zero real constant and $\hat g=\frac{2}{3\la \sin(\frac{x}{3})^2} g_{\R}$ be the conformally rescaled metric of $g_{\R}$ in (\ref{m1f}).
Then $\text{Ric}(\hat g)=\la \hat g$. In particular, $\hat g$ is an Einstein metric with constant scalar curvature $R=5\la$. 
\end{proposition}

When $c^2=9$, we take $h(x)=\cos(3x)$.
This gives the real Pfaffian system $\cP_{\R}$ spanned by the 1-forms
\begin{align*}
o_1=\si_1-\cos(3x)\der q,\qquad
o_2=\si_2-\der x,\qquad
o_3=\si_3+3\sin(3x)\der q,
\end{align*}
with $o_4=-\chi_3=-\der x$ and $o_5=\chi_4=\cos(3x)\der q$.

This Pfaffian system annihilates the vector fields 
\begin{align*}
X^1=\frac{\partial}{\partial x}+E^2,\qquad
X^2=\frac{\partial}{\partial q}+\cos\left(3x\right)E^1-3\sin\left(3x\right)E^3,
\end{align*}
which span the real distribution $\cD_{\R}$.
Reversing signs in $o_1$ and $o_2$, we take $\tilde \cP_{\R}$ to be the Pfaffian system given by the 1-forms
\begin{align*}
\tilde o_1=\si_1+\cos(3x)\der q,\qquad
\tilde o_2=\si_2+\der x,\qquad
\tilde o_3=o_3,
\end{align*}
which annihilate the ``sign-reversed" distribution $\tilde \cD_{\R}$.
The associated real metric given by 
\begin{align*}
g_{\R}&=-2(o_4 o_4+o_5o_5)+\frac{1}{2}\tilde{o}_1\tilde{o}_1+\frac{1}{2}\tilde{o}_2\tilde{o}_2-\frac{1}{6}\tilde{o}_3 \tilde{o}_3
\end{align*}
is conformally flat.
We have the following theorem:
\begin{theorem}\label{real3a}
A representative metric of Nurowski's conformal class associated to the maximally symmetric real distribution $\cD_{\R}$ given by the span of 
\begin{align*}
X^1=\frac{\partial}{\partial x}+E^2,\qquad
 X^2=\frac{\partial}{\partial q}+\cos\left(3x\right)E^1-3\sin\left(3 x\right)E^3,
\end{align*}
is given by 
\begin{align}\label{m2f}
g_{\R}&=-2(o_4 o_4+o_5o_5)+\frac{1}{2}\tilde{o}_1\tilde{o}_1+\frac{1}{2}\tilde{o}_2\tilde{o}_2-\frac{1}{6}\tilde{o}_3 \tilde{o}_3\\
&=\chi_1^2+\chi_2^2-(\chi_3^2+\chi_4^2)-\frac{1}{2}(\om_1^2+\om_2^2)-\frac{1}{6}\om_3^2,\nonumber
\end{align}
where
$\tilde o_1$, $\tilde o_2$, $\tilde o_3$ are 1-forms in the ``sign-reversed" Pfaffian system $\tilde \cP_{\R}$, which span the complementary bundle $\tilde \cD_{\R}^{\perp}$. The conformal structure $[g_{\R}]$ is conformally flat and the distribution $\cD_{\R}$ has the split real form of $G_2$ as its group of symmetries.
\end{theorem}

Likewise for the metric obtained in Theorem \ref{real3a}, we are able to find an Einstein representative in the conformal class. It is given by the following. 
\begin{proposition}\label{einb}
Let $\la$ be a non-zero real constant and $\hat g=-\frac{6}{\la \cos(\theta)^2} g_{\R}$ be the conformally rescaled metric of $g_{\R}$ in (\ref{m2f}).
Then $\text{Ric}(\hat g)=\la \hat g$. In particular, $\hat g$ is an Einstein metric with constant scalar curvature $R=5\la$. 
\end{proposition}

 \section{Complexified $SU(2)$ Pfaffian systems}\label{max4}
 In this section, we consider the half-complexified setting where $f(\theta)=-\sin(\theta)$ and $h(x)=e^{c i x}$ with $c^2=\frac{1}{9}$ or $9$. In other words we fix the real left-invariant vector fields on $SU(2)$ and allow $h(x)$ to be complex. This example was considered first in earlier versions of the manuscript because it is a homogeneous example of a $(2,3,5)$-distribution. Completing the 1-forms ($o_1$,$o_2$,$o_3$,$o_4$,$o_5$) to a coframing for Nurowski's metric gives only constants appearing in the structure equations and therefore defines a homogeneous structure. However, following the peer review process, we have decided to focus on the rolling distribution and $SU(2)$ Pfaffian systems in the real setting, presenting it first and we will only give brief mention to the complexified results obtained here. The results are basically similar to the previous section, and we highlight the result about Ricci-flatness for the metric considered here upon conformally rescaling the metric $\tilde g$.

When $c=\mp \frac{1}{3}$, taking $f(\theta)=-\sin(\theta)$ and $h(x)=e^{c i x}$ gives us complex distributions $\cD^{\mp}$ spanned by the vector fields
\begin{align*}
X^1=\frac{\partial}{\partial x}+E^2,\qquad X^{2{\mp}}=\frac{\partial}{\partial q}+ e^{\mp \frac{i}{3}x} E^1\mp \frac{i}{3} e^{\mp \frac{i}{3}x} E^3.
\end{align*}
We have the following result. 
\begin{theorem}\label{c1m}
A representative metric of Nurowski's conformal class associated to the maximally symmetric complex distributions $\cD^{\mp}$ given by the span of 
\begin{align*}
X^1=\frac{\partial}{\partial x}+E^2,\qquad X^{2{\mp}}=\frac{\partial}{\partial q}+ e^{\mp \frac{i}{3}x} E^1\mp \frac{i}{3} e^{\mp \frac{i}{3}x} E^3,
\end{align*}
is given by 
\begin{align*}
g_{\mp}&=-2(o_4 o_4+o_5o_5)+\frac{1}{2}\tilde{o}_1\tilde{o}_1+\frac{1}{2}\tilde{o}_2\tilde{o}_2+\frac{1}{2}\tilde{o}_3 \tilde{o}_3+(o_1o_1+o_2o_2+o_3o_3)\\
&=-(\om_4 \om_4+\om_5 \om_5)+(\si_1\si_1+\si_2\si_2)+\frac{1}{2}(\om_1\om_1+\om_2\om_2)+\frac{3}{2}\om_3\om_3.
\end{align*}
Here
\begin{align*}
o_1&=\si_1-e^{\mp\frac{i}{3}x}\der q,\qquad
o_2=\si_2-\der x,\qquad
o_3=\si_3\pm \frac{i}{3}e^{\mp \frac{i}{3}x}\der q,
\end{align*}
are 1-forms in the Pfaffian system $\cP^{\mp}$, which span $(\cD^{\mp})^{\perp}$ and 
\begin{align*}
\tilde{o}_1&=\si_1+e^{\mp\frac{i}{3}x}\der q, \qquad
\tilde{o}_2=\si_2+\der x,\qquad
\tilde{o}_3=o_3,
\end{align*}
are 1-forms in the Pfaffian system $\tilde \cP^{\mp}$, which span $(\tilde \cD^{\mp})^{\perp}$. We also have 
\begin{align*}
o_4&=-\der x,\qquad o_5=e^{\mp\frac{i}{3}x}\der q.
\end{align*}
The metric $g_{\mp}$ is conformally flat.
\end{theorem}

We can further investigate the question of whether the conformal structure $[g_{\mp}]$ associated to $\cD^{\mp}$ admits a Ricci-flat representative by solving for the equation 
\[
\text{Ric}(\Om(x)^2 g_{\mp})=0. 
\] 
When $h(x)=e^{-\frac{i}{3}x}$, we find that the conformal factor $\Om(x)$ satisfies the complex first order differential equation 
\[
\Om_x=\frac{i}{3}\Om, 
\]
which has complex solutions given by $\Om(x)=\ka e^{\frac{i}{3}x}$.
When $h(x)=e^{\frac{i}{3}x}$, the conformal factor $\Om(x)$ satisfies the complex differential equation 
\[
\Om_x=-\frac{i}{3}\Om, 
\]
which has complex solutions given by $\Om(x)=\ka e^{-\frac{i}{3}x}$. The appearance of the complex factor in the conformal rescaling suggests that the conformal factor can be interpreted as a phase factor. 

\begin{proposition}
The flat conformal structure $[g_{\mp}]$ determined by $g_{\mp}$ in Theorem \ref{c1m} is Ricci-flat when multiplied by the square of the phase factor $\Om(x)=\ka e^{{\pm}\frac{i}{3}x}$ (the sign chosen corresponding to the sign of $\cP^{\mp}$), giving us a Ricci-flat representative of $[g_{\mp}]$. 
\end{proposition}

We now investigate the other case where $b=1$ and $c^2=9$, again following the preceding section. We have the following result.
\begin{theorem}\label{c2ma}
A representative metric of Nurowski's conformal class associated to the maximally symmetric complex distribution $\cD^{\mp}$ given by the span of 
\begin{align*}
X^1=\frac{\partial}{\partial x}+E^2,\qquad X^{2{\mp}}=\frac{\partial}{\partial q}+e^{\mp 3 i x}E^1\mp 3 i e^{\mp 3 i x} E^3,
\end{align*}
is given by 
\begin{align*}
g_{\mp}&=-2(o_4 o_4+o_5o_5)+\frac{1}{2}\tilde{o}_1\tilde{o}_1+\frac{1}{2}\tilde{o}_2\tilde{o}_2-\frac{1}{6}\tilde{o}_3 \tilde{o}_3\\
&=\si_1\si_1+\si_2\si_2-o_4o_4-o_5o_5-\frac{1}{2}(o_1o_1+o_2o_2)-\frac{1}{6}o_3o_3.
\end{align*}
Here
\begin{align*}
o_1=\si_1-e^{\mp3i x}\der q,\qquad
o_2=\si_2-\der x,\qquad
o_3=\si_3\pm3i e^{\mp3ix}\der q
\end{align*}
are 1-forms in the Pfaffian system $\cP^{\mp}$ and
\begin{align*}
\tilde o_1=\si_1+e^{\mp3i x}\der q,\qquad
\tilde o_2=\si_2+\der x,\qquad
\tilde o_3=o_3
\end{align*}
are 1-forms in the ``sign-reversed" Pfaffian system $\tilde \cP^{\mp}$, which span the complementary bundle $(\tilde \cD^{\mp})^{\perp}=(TM/\tilde \cD^{\mp})^*$.
 We also have 
\begin{align*}
o_4&=-\der x,\qquad o_5=e^{\mp3 ix}\der q.
\end{align*}
The metric $g_{\mp}$ is conformally flat.
\end{theorem}

Interestingly, we are unable to find a Ricci-flat representative in the conformal class of the metric in Theorem \ref{c2ma}. Instead, Proposition \ref{einb} holds and when $\la=-6$, we have the following:
\begin{proposition}
Let $\hat g=\frac{1}{\cos(\theta)^2}g_{\mp}$ be the conformally rescaled metric of $g_{\mp}$ in Theorem \ref{c2ma}.
Then $\text{Ric}(\hat g)=-6\hat g$. In particular, $\hat g$ is an Einstein metric with constant scalar curvature $R=-30$. 
\end{proposition}

The asymmetry in the choices of $f(\theta)$ and $h(x)$ seems to lead to the asymmetry in the conformal rescaling of the metrics in Theorems \ref{c1m} and \ref{c2ma}. Indeed it is unclear whether a Ricci-flat representative metric for $g_{\mp}$ in Theorem \ref{c2ma} can be found.
 
To finish the discussion, in the work of \cite{R21b}, the parametrisation taking $f(\theta)=e^{i \theta}$ and $h(x)=e^{i c x}$ where $c=\pm \frac{1}{3}$ or $\pm 3$ was considered. In this fully complexified setting, the author was able to find the change of coordinates to the flat Cartan distribution and write down the vector fields that bracket-generate split $\frak{g}_2$ (see \cite{R21b}).  Indeed, we can think of this setting as complexifying the change of coordinates we obtained earlier in Section \ref{max2}, replacing $(\theta, x)$ with $(i \theta, i x)$. The complex vector fields satisfying the $su(2)$ commutation relations that we obtained however are no longer the real left-invariant vector fields on $SU(2)$. 

\section{Doubling the $SU(2)$ symmetry}\label{max5}

Under symmetry considerations, we can also generalise the results in Section \ref{max3} on $SU(2)$ Pfaffian systems by supposing that the sphere element $o_4o_4+o_5o_5$ arises from another copy of $SU(2)$. Let $(\psi,\theta,\phi,v,x,q,z)$ be coordinates parametrising $S^3 \times S^3 \times \R$, with the left-invariant 1-forms on the second copy of $S^3$ given by
\begin{align*}
\si_4&=\sin(v) \der x-\cos(v) \sin(x) \der q,\\
\si_5&=\cos(v)\der x+\sin(v) \sin(x) \der q,\\
\si_6&=-\der v-\cos(x) \der q. 
\end{align*}
We shall take
\begin{align*}
\tau_1=\si_2+i \si_1,\qquad \tau_2=\si_2-i \si_1,\qquad \tau_4=\si_5+i \si_4,\qquad \tau_5=\si_5-i \si_4.
\end{align*}
We have
\begin{align*}
\tau_1&=e^{i \psi}(\der \theta-i\sin(\theta) \der \phi),\qquad \tau_2=e^{-i \psi}(\der \theta+i\sin(\theta) \der \phi),\\
\tau_4&=e^{i v}(\der x-i\sin(x) \der q),\qquad \tau_5=e^{-i v}(\der x+i\sin(x) \der q).
\end{align*}
Consider the rank $3$ Pfaffian system on $S^3 \times S^3 \times \R$ given by
\begin{align*}
\eta_1&=\si_2+i \si_1+\al e^{-i(v-\psi+z)}(\si_5+i \si_4),\\
\eta_2=\bar{\eta}_1&=\si_2-i \si_1+\al e^{i(v-\psi+z)}(\si_5-i \si_4),\\
\eta_3&=-\der (v-\psi+z)+\si_3-\si_6.
\end{align*}
The 1-forms annihilate the vector fields spanned by
\begin{align*}
\partial_z-\al(E^1-iE^2)+\ga E^3+e^{i(\psi-z-v)}(E^4-iE^5)+(\ga+1)E^6
\end{align*} and its complex conjuate. Here $\ga$ is a constant. 
Let $M^5\subset S^3\times S^3\times \R$ be the submanifold given by $\{z=v=\psi\}$, and the 1-forms when restricted to $M^5$ give
\begin{align}\label{pf1}
\eta_1&=\si_2+i \si_1+\al e^{-iz}(\si_5+i \si_4)=\si_2+i\si_1+\al(\der x-i \sin(x) \der q),\nonumber\\
\eta_2=\bar{\eta}_1&=\si_2-i \si_1+\al e^{iz}(\si_5-i \si_4)=\si_2-i\si_1+\al(\der x+i \sin(x) \der q),\\
\eta_3&=-\der z+\si_3-\si_6=-\der z-\cos(\theta) \der \phi+\cos(x) \der q,\nonumber
\end{align}
where we have used the fact that on $M^5$, 
\begin{align*}
\si_4&=\sin(z) \der x-\cos(z) \sin(x) \der q,\\
\si_5&=\cos(z)\der x+\sin(z) \sin(x) \der q,\\
\si_6&=-\der z-\cos(x) \der q. 
\end{align*}
We have the following theorem:
\begin{theorem}
On $M^5$, the $(2,3,5)$-distribution associated to the rank 3 Pfaffian system spanned by the 1-forms $\{\eta_1, \eta_2, \eta_3\}$ in (\ref{pf1}) is maximally symmetric when $\al^2=\frac{1}{9}$ or $\al^2=9$, i.e.\ when $\al=\pm\frac{1}{3}$ or $\pm 3$. In the case when $\al=\pm\frac{1}{3}$, a representative metric in Nurowski's conformal structure associated to this distribution is given by
\begin{align}\label{metric1}
g=9 \eta_1\eta_2+3\eta_3\eta_3-18\tau_1\tau_2+2\tau_4\tau_5.
\end{align}
 In the case when $\al=\pm 3$, a representative metric in Nurowski's conformal structure associated to this distribution is given by
\begin{align*}
g=\eta_1\eta_2+3\eta_3\eta_3+2\tau_1\tau_2-18\tau_4\tau_5.
\end{align*}
The two metrics are equivalent via the substitution
$(z,\theta,\phi,x,q)\mapsto(-z,x,q,\theta,\phi)$. Both are conformally flat with vanishing Weyl tensor. Under a conformal rescaling, 
the first (resp. second) metric
\begin{align*}
\hat g=\frac{4}{\cos(\theta)^2}g\qquad \left(resp.~ \hat g=\frac{4}{\cos(x)^2}g\right)
\end{align*}
is Einstein with 
\begin{align*}
\text{Ric}\left(\hat g\right)=\frac{1}{12}\hat g,\qquad \left(resp.~ \text{Ric}\left(\hat g\right)=\frac{1}{12}\hat g \right).
\end{align*}
\end{theorem}
From this form of the metric, we can readily recover the 1-forms of the $SU(2)$ Pfaffian system in Section \ref{max3}. The distribution spanned by the kernel of the 1-forms also annihilate the combinations
\begin{align*}
\frac{1}{2i}\left(\eta_1-\bar{\eta}_1\right)=\si_1-\al \sin(x) \der q,\qquad \frac{1}{2}\left(\eta_1+\bar{\eta}_1\right)=\si_2+\al \der x,\qquad \eta_3=\si_3+\cos(x) \der q.
\end{align*}
This Pfaffian system is invariant under the substitution $(x,q,\al)\mapsto (-x, q,-\al)$.

To relate the distribution back to the rolling surfaces distribution of An and Nurowski, observe that the elements given by $18\tau_1\tau_2$ and $2\tau_4\tau_5$ in (\ref{metric1}) describe the metrics on a pair of spheres. The element with coefficient $2$ corresponds to the sphere of radius $\sqrt{2}$ while the element with coefficient $18$ corresponds to the sphere of radius $3\sqrt{2}$. The ratio of their curvatures are in the form $9:1$.

Replacing $z$ by $\psi$, we shall let 
\begin{align*}
E^4&=\sin(\psi)\partial_x+\cos(\psi)\left(\cot(x)\partial_{\psi}-\frac{1}{\sin(x)}\partial_q\right),\\
E^5&=\cos(\psi)\partial_x-\sin(\psi)\left(\cot(x)\partial_{\psi}-\frac{1}{\sin(x)}\partial_q\right).
\end{align*}
The vector fields $E^3$, $E^4$, $E^5$ form another set of left-invariant vector fields of $SU(2)$. 
The vector fields annihilating the 1-forms (\ref{pf1}) are given by
\begin{align}\label{v1a}
X^1&=2\cot(x)E^3-\al(E^1-iE^2)+e^{-i\psi}(E^4-iE^5),\\
X^2=\bar{X}^1&=2\cot(x)E^3-\al(E^1+iE^2)+e^{i\psi}(E^4+iE^5).\label{v1b}
\end{align}
Using the fact that
\begin{align*}
\partial_x=\sin(\psi)E^4+\cos(\psi)E^5 \text{~and~} -\frac{1}{\sin(x)}\partial_q-\cot(x)E^3=\cos(\psi)E^4-\sin(\psi)E^5,
\end{align*}
we find that $-\frac{1}{2i}(X^1-X^2)$ gives the vector field
\begin{align*}
\partial_x-\al E^2, 
\end{align*}
while $-\frac{1}{2}\sin(x)(X^1+X^2)$ gives the vector field
\begin{align*}
\partial_q+\al \sin(x) E^1-\cos(x) E^3. 
\end{align*}
This pair of vector fields describe a $su(2)$-valued connection over the sphere with coordinates $(x,q)$. 

We can interpret the complex vector fields $X^1$ and $X^2$ that span the distribution $\cD$ as determining a non-Abelian $SU(2)$ spin connection on a rank 2 spinor bundle $\cE'$ over $\Si'^2$. Using the $U(1)$ symmetry of the circle twistor bundle over the product of these two surfaces $\Si^2 \times \Si'^2$ (see \cite{AN14}), we show that the distribution can be rotated to define a spin connection on $\cE$ over $\Si^2$ with $\Si \hookrightarrow SU(2) \subset\cE$. 

We show the following:
\begin{theorem}
By passing to the An-Nurowski circle twistor distribution over the product of two spheres $\Si ^2\times \Si'^2$, the $su(2)$-valued connection over the sphere 
$\Si'^2$ with coordinates $(x,q)$ given by
\begin{align*}
\partial_x-\al E^2,\qquad 
\partial_q+\al \sin(x) E^1-\cos(x) E^3,
\end{align*}
can be rotated to a $su(2)$-valued connection over the other sphere $\Si^2$ with coordinates $(\theta, \phi)$ and components given by
\begin{align*}
\partial_\theta-\frac{1}{\al} E^5,\qquad 
\partial_\phi+\frac{1}{\al} \sin(\theta) E^4-\cos(\theta) E^3.
\end{align*}
\end{theorem}

\begin{proof}
The rotated connection can be obtained in the following way. Let $X^1$ and $X^2$ be the vector fields (\ref{v1a}) and (\ref{v1b}), which define the connection over $\Si'^2$. Take the combinations
\begin{align*}
Y^1=\frac{1}{2}(e^{i\psi}X^1+e^{-i\psi}X^2),\quad Y^2=\frac{1}{2i}(e^{i\psi}X^1-e^{-i\psi}X^2).
\end{align*}
Then we have the vector fields
\begin{align*}
-\sin(2\psi)Y^1+\cos(2\psi) Y^2={}&\al \partial_\theta+\cot(x)\sin(\psi)\partial_{\psi}-\cos(\psi)\partial_x+\frac{\sin(\psi)}{\sin(x)}\partial_q,\\
\cos(2\psi)Y^1+\sin(2\psi) Y^2={}&\frac{\al}{\sin(\theta)} \partial_\phi-(\cot(x)\cos(\psi)+\al\cot(\theta))\partial_{\psi}\\
&-\sin(\psi)\partial_x-\frac{\cos(\psi)}{\sin(x)}\partial_q.
\end{align*}
Now let $\tilde \psi=-\psi$. This gives $\partial_{\tilde \psi}=-\partial_{\psi}$ and we have
\begin{align*}
\sin(2\tilde \psi)Y^1+\cos(2\tilde \psi) Y^2={}&\al \partial_\theta+\cot(x)\sin(\tilde \psi)\partial_{\tilde \psi}-\cos(\tilde \psi)\partial_x-\frac{\sin(\tilde \psi)}{\sin(x)}\partial_q,\\
={}&\al \partial_\theta-E^5,\\
\cos(2\tilde \psi)Y^1-\sin(2\tilde \psi) Y^2={}&\frac{\al}{\sin(\theta)} \partial_\phi+(\cot(x)\cos(\tilde \psi)+\al\cot(\theta))\partial_{\tilde \psi}\\
&+\sin(\tilde \psi)\partial_x-\frac{\cos(\tilde \psi)}{\sin(x)}\partial_q,\\
={}&\frac{\al}{\sin(\theta)} \partial_\phi+E^4+\al\cot(\theta)\partial_{\tilde \psi}\\
={}&\frac{\al}{\sin(\theta)} \partial_\phi+E^4-\al\cot(\theta)E^3.
\end{align*}
This combination gives, up to proportionality, the vector fields
\begin{align*}
\partial_\theta-\frac{1}{\al} E^5\quad\mbox{and}\quad 
\partial_\phi+\frac{1}{\al} \sin(\theta) E^4-\cos(\theta) E^3,
\end{align*}
which define the connection over $\Si^2$.

\end{proof}

\end{document}